\documentclass{article}
\usepackage{amsmath}
\usepackage{float}
\usepackage{graphicx}
\graphicspath{ {imgs/} }
\usepackage{amsthm}
\usepackage{amsfonts}
\begin{document}

\title{Pair Correlations in Uniform Countable Sets}
\author{C. Schildkraut, S. Raman}
\maketitle
\newtheorem{theorem}{Theorem}[section]
\newtheorem{corollary}[theorem]{Corollary}
\newtheorem{lemma}{Lemma}

\theoremstyle{definition}
\newtheorem{definition}{Definition}[section]

\newcommand{\Beta}{\mathrm{B}}
\def\Tau{\mathop{\mathrm{T}}}
\newcommand{\zerovector}{\mathbf{0}}
\newcommand{\Vu}{\mathbf{u}}
\newcommand{\Vv}{\mathbf{v}}
\newcommand{\Vx}{\mathbf{x}}
\newcommand{\Vy}{\mathbf{y}}
\newcommand{\fT}{\mathcal{T}}
\newcommand{\fS}{\mathcal{S}}
\newcommand{\fU}{\mathcal{U}}
\newcommand{\fZ}{\mathcal{Z}}
\bibliographystyle{plain}

\begin{abstract}
We determine the pair correlations of countable sets $T \subset \mathbb{R}^n$ satisfying natural equidistribution conditions. The pair correlations are computed as the volume of a certain region in $\mathbb{R}^{2n}$, which can be expressed in terms of the incomplete Beta function.  For $n=2$ and $n=3$ we give closed form expressions, and we obtain an expression in the $n \rightarrow \infty$ limit. We also verify that sets of lattice points and primitive lattice points satisfy the required distribution criteria. 
\end{abstract}

\section{Introduction}

In many areas of geometry it is important to determine the spatial statistics of various discrete sets of points in $n$-dimensions. We study an interesting property of the statistics of integer lattice points $\mathbb{Z}^{n}$ (and subsets of the lattice points), demonstrating an interesting relationship between probability and number theory. In particular, we calculate the \textit{pair-correlations} between the points of a uniform countable set. 

Given a countable set of points $T \subset \mathbb{R}^n$ satisfying certain natural uniform distribution conditions (as are satisfied by many standard discrete sets, like the lattice points and the primitive lattice points), we consider the subset of points $T$ whose magnitudes are $\leq r$. Pick two points in this set $T \cap B_n(r)$ (where $B_n(r)$ is the $n$-ball centered at $\zerovector$ of radius $r$), and let the distance between them be $\lambda R$. We want to find the probability density function of $\lambda$ in the limit as $R \to \infty$. To do this, we first prove the heuristic relating the number of points in $T\times T$ within a region in $\mathbb{R}^{2n}$ to its volume. After this, we make use of our volume heuristic to reduce the computation of the required probability density function to a volume calculation. 

\subsection{Conditions for Equidistribution}

For our result to hold good, we demand that T satisfy the following equidistribution conditions: 

\vspace{2mm}
\noindent \textbf{I. Angular Condition.}
We require that $T$ becomes uniformly distributed on the sphere $\mathbb S^{n-1} = \partial B_n (1)$ in the sense that, for any continuous function $f$ on $\mathbb S^{n-1}$ extended to the whole of $\mathbb{R}^n$ by $f(\Vx) \equiv f \left( \hat{\Vx} \right) $ (where $\hat{\Vx}$ is the unit vector pointing in the direction of $\Vx$), we have  

\begin{equation}\label{AngularCriterion}
 \lim_{r \rightarrow \infty}  \frac{1}{\left| T \cap B_n(r) \right|} \sum_{\Vx \in T \cap B_n(r)} f(\Vx) = \int_{\mathbb S^{n-1}} f(\Vx)d\mu_P(\Vx),
\end{equation}

\noindent where $\mu_P$ denotes the Lebesgue probability measure on $\mathbb S^{n-1}$.

\vspace{2mm}
\noindent \textbf{II. Radial (Growth) Condition.}
We moreover require that $T$ satisfies the following (natural) radial counting: letting $N(r) = N(T, r)$ be the number of points $\Vx$ in $T$ such that $|\Vx| < r$, we take 

\begin{equation} \label{RadialCriterion}
\frac{N(r\lambda)}{N(r)} \sim \lambda^n .
\end{equation}

\noindent We wish to show that, for some $\Psi \subset \mathbb{R}^{2n}$ with $\Psi(r)$ being $\Psi $ dilated by $r$ about $\zerovector$, $|(T \times T) \cap \Psi(r) | $ grows like $N(r)^2$, in that the limit as $r\to\infty$ of their quotient tends to $m(\Psi)/m(B_n)^2$, where $m$ is the (unnormalized) Lebesgue measure in Euclidean space.  

\section{The Volume Heuristic for Equidistributed Sets}

To calculate pair correlations of sets $T$ satisfying these equidistribution conditions (\ref{AngularCriterion}) and (\ref{RadialCriterion}), we require a heuristic relating the number of ordered pairs of points of $T$ in a certain region in $\mathbb{R}^{2n}$ to the volume of that region. 

We approach the proof of this in three steps. First, we provide a much simpler equidistribution criterion equivalent to the original integral one. Then, we apply this new criterion to prove an intermediate result relating the number of points of $T$ within regions in $\mathbb{R}^n$. Finally, we use an argument involving Cartesian products to extend this to $\mathbb{R}^{2n}$. 

First, we cite a theorem which states the equivalence of the angular equidistribution criterion to a much simpler statement. This classical result roughly asserts that a set is equidistributed in the angular sense if it looks equally dense in any direction.  

\begin{theorem}[Kuipers and Niederreiter 1974]

Define 
\begin{equation}
\label{WedgeDef}
W(X,r) = \left\{\Vx\ \big|\ |\Vx|<r, \frac{\Vx}{|\Vx|} \in X\right\}  .
\end{equation}
Given a countable set $T \subset \mathbb{R}^n$, the above equidistribution condition is true iff the following holds: A set $T$ is equidistributed over $\mathbb{S}^{n-1}$ if for all $X \subseteq \mathbb{S}^{n-1} $, we have
\begin{equation} 
\label{KuipersSimplification}
 \lim_{r \to \infty} \frac{ | T \cap W(X, r) |  }{N(r)} = \mu_P(X) = \frac{ m( W(X, 1) )}{m(B_n) }, 
\end{equation}

\noindent where $\mu_P$ is the probability measure on $\mathbb{S}^{n-1}$.

\end{theorem}

\noindent This theorem is proven in \cite{Kuipers1974}. Making use of this result, we now prove an preliminary theorem from which the main volume heuristic can be derived. Loosely speaking, this theorem states that any sufficiently well-behaved region in $\mathbb{R}^n$ contains roughly as many points of a given equidistributed set as the numerical value of its volume (in an appropriate limit).

\begin{theorem} 

For a given compact set $\Omega \subset \mathbb{R}^n$ such that $\forall \ \Vx \in \Omega$, $\sigma\Vx \in \Omega\ \forall\ 0\leq \sigma\leq 1$ and that $m(\partial \Omega) = 0$, let $\Omega(r)$ be $\Omega$ dilated by $r$ about $\zerovector$. Also, let $T \subset \mathbb{R}^n$ satisfy the conditions given at the beginning of the paper. We then have
\begin{equation}
\label{IntermediateEquidistributionResult}
\frac{\left| T \cap \Omega(r) \right|}{N(r)} \sim \frac{m(\Omega)}{m(B_{n})},
\end{equation} 
where $N(r)$ is as defined at the beginning of the paper. 

\begin{proof}

Essentially, our proof involves breaking up $\Omega$ directionally into small cones, applying the assumed angular equidistribution condition (Theorem 2.1), and summing the results. For a given unit $n$-vector $\mathbf{v}$, let $\lambda({\mathbf{v}}) = \sup_{ \lambda \Vv \in \Omega} (\lambda) $. Note that $\sigma\mathbf{v} \in \Omega\ \forall\ 0\leq \sigma < \lambda(\mathbf{v})$. 

\vspace{2mm}

Define the diameter of a compact Borel set $X \subseteq \partial B_n$ to be

\begin{equation} \label{DiameterDef}
\Delta(X) = \sup_{\Vx, \Vy \in X}  \Vert \Vx - \Vy \Vert  .
\end{equation} 

\noindent In other words, the diameter is the least upper bound on the distance between any two points in $X$. Consider a partition $\mathcal{P}(\delta)$ of $\partial B_n$ into Borel sets $X$ with $\Delta(X)<\delta \ \forall\ X \in \mathcal{P}(\delta)$ for some $\delta > 0$. We can approximate our region $\Omega$ by the union of all of the $W\left(X,\lambda(\mathbf{x})\right)$ for some $\mathbf{x} \in X$ for each $X \in \mathcal{P}(\delta)$. Note that these regions are disjoint (because $\mathcal{P}(\delta)$ is by definition a disjoint union of sets). We see that as $\delta \rightarrow 0$, the union of these regions $W\left(X,\lambda(\mathbf{x})\right)$ tends towards $\Omega$ itself. As such, let $\Omega_{\delta}$ be the approximation of $\Omega$ by a partition $\mathcal{P}(\delta)$; then, $\lim_{\delta \to 0} \Omega_{\delta}(r) = \Omega(r)$. We have that

\begin{equation}
\label{SumOverPartition}
\left|T \cap \Omega_\delta(r) \right| =  \bigg( \sum_{X \in \mathcal{P}(\delta)} \left|T \cap W(X,r\lambda(\mathbf{x}))\right| \bigg) \ ;
\end{equation}

\noindent dividing by $N(r)$, taking $r \to \infty$, and then $\delta \to 0$, we see that 

\begin{equation} \label{BeforeRadial}
\lim_{r \to \infty}  \frac{ \left|T \cap \Omega(r) \right|}{N(r)} = \lim_{\delta \rightarrow 0} \lim_{r \to \infty} \frac{1}{N(r)}  \sum_{X \in \mathcal{P}(\delta)} \left|T \cap W(X,r\lambda(\mathbf{x}))\right| \ . 
\end{equation}
 
\noindent Also, condition (\ref{AngularCriterion}) (equivalently (\ref{KuipersSimplification})) gives

$$\lim_{r\to\infty} \frac{\left| T \cap W(X,r\lambda)\right|}{N(r\lambda)} = \frac{m(W(X,1))}{m(B_{n})}$$

\noindent for all $X\in \mathcal{P}(\delta)$ and any positive $\lambda$. Using condition (\ref{RadialCriterion}), we obtain

\begin{equation} \label{RadialResult}
\lim_{r\to\infty} \frac{\left|T \cap W(X,r\lambda)\right|}{N(r)} = \frac{m(W(X,\lambda))}{m(B_{n})}.
\end{equation}

\noindent Now, we see that the asymptotic expression can be used as $r$ grows large; this enables us to write

$$\lim_{r\to\infty} \frac{\left| T \cap \Omega(r) \right|}{N(r)} = \lim_{\delta\to 0}\sum_{X \in \mathcal{P}(\delta)} \frac{m(W(X,\lambda(\mathbf{x}) ))}{m(B_{n})} $$

$$\lim_{r\to\infty} \frac{\left| T \cap \Omega(r) \right|}{N(r)} = \frac{m(\Omega)}{m(B_{n})},$$

\noindent which is exactly (\ref{IntermediateEquidistributionResult}).
\end{proof}

\end{theorem}

We now use the above result to prove the main theorem, a result showing the main volume heuristic. This statement naturally restates the previous result (Theorem 2.2) to $2n$-dimensions, but the details of this extension are quite complicated. It turns out that $\mathbb{R}^{2n}$ is far more convenient in which to work with pair correlations (for each distance, there are two independent points with $n$ degrees of freedom), which is why the following theorem is of great importance.

\begin{theorem}[Volume Heuristic]

Let $K$ be the region in $2n$-space (with coordinates $x_1,x_2,...,x_{2n}$) defined by

\begin{equation} \label{Bn2Definition}
K = \left\{ \Vx = (x_1, \cdots, x_{2n}) \ \ \bigg| \ \sum_{i=1}^n x_i^2 \leq 1\mathrm{\ and\ }\sum_{i=1}^n x_{n+i}^2 \leq 1 \right\},
\end{equation}

\noindent and let $\Psi \subseteq K$ be a compact region. Furthermore, let $K(r)$ be $K$ dilated by $r$, and let $\Psi(r)$ be $\Psi$ dilated by $r$. Also, define $T$ as before, and let $T \times T$ be the set of vectors formed by the concatenation of the two $n$-vectors $(x_1, \dots x_n)$ and $(x_{n+1}, \dots x_{2n})$, each of which is in $T$. Then, we have that

\begin{equation}
\label{VolumeHeuristicStatement}
\lim_{r\to\infty} \frac{\big|(T \times T) \cap \Psi(r) |}{\big| (T \times T) \cap K(r) \big|} = \frac{m(\Psi)}{m(B_n)^2}.
\end{equation}

\begin{proof}

Our strategy here is slightly different. Since we are trying to extract a $2n$-dimensional result from criteria in $n$ dimensions, we partition $\mathbb{R}^{2n}$ into cubes, apply Theorem 2.2 to the $n$-dimensional slices of these cubes, and take the Cartesian product. Before this, however, we prove an important lemma which shows that the contribution of points on the boundaries of these cubes is asymptotically negligible. 

Let $\delta$ be a real number, and $\mathcal{H_\delta}$ be the following partition of $\mathbb{R}^{2n}$ into hypercubes of side length $\delta$:

\begin{equation} \label{HDefinition}
\mathcal{H}_\delta = \left( \bigcup_{i_1 = -\infty}^{\infty} \left( \delta  i_1 - \delta, \delta  i_1  \right] \right) \times \cdots \times \left( \bigcup_{i_{2n} = -\infty}^{\infty} \left( \delta  i_{2n} - \delta , \delta  i_{2n} \right] \right),
\end{equation}

\noindent where $\times $ is the Cartesian product. In other words, a part of $\mathcal{H}_{\delta}$ would be $\left( \delta  i_1 - \delta , \delta i_1  \right] \times \cdots  \left( \delta  i_{2n} - \delta , \delta  i_{2n}  \right] $ for some integer $2n$-vector $(i_1, \dots i_{2n})$. Let the set of all such hypercubes contained entirely in $\Psi$ be $\mathfrak{C}_{\delta}$. We first prove the following lemma: 

\begin{lemma}[Boundary Points Lemma]

Consider a face of an $n$-box $C$ with sides parallel to the axes; call it $\mathcal{F}$, and let $\mathcal{F}(r)$ be $\mathcal{F}$ dilated by $r$. Then $| T \cap \mathcal{F}(r) | = o(N(r)) $.  

\begin{proof} 

Without loss of generality, assume that $x_1$ is held constant for this face and all other coordinates may vary within the hypercube.  Consider the region $\Gamma(r)$ bounded by the lines from the origin to the vertices of $\mathcal{F}(r)$ and bounded face $\mathcal{F}(r)$ itself (we let $\Gamma$ be $\Gamma(r)$ dilated by $1/r$ and define $\mathcal{F}$ similarly). Note that we may apply Theorem 2.2 to $\Gamma(r)$, which can easily be seen. With this in mind, dilate $\Gamma(r)$ by $1 - \epsilon$ for some $\epsilon > 0$, and apply Theorem 2.2 to obtain

\begin{equation} \label{OneMinusEpsilon}
\lim_{r \to \infty} |T \cap \Gamma(r(1 - \epsilon))| = \lim_{r \to \infty} \frac{m(\Gamma) N(r(1 - \epsilon)) }{m(B_n)}.
\end{equation}

\noindent Also, dilate $\Gamma(r)$ by $1+ \epsilon$ to get $\Gamma(r(1+\epsilon))$; then, we also obtain 

\begin{equation} \label{OnePlusEpsilon}
\lim_{r \to \infty} |T \cap \Gamma(r(1+\epsilon))| = \lim_{r \to \infty} \frac{m(\Gamma)N(r(1 + \epsilon))}{m(B_n)}.
\end{equation}

\noindent Now, $|T \cap \mathcal{F}(r)| \leq |T \cap \Gamma(r(1 + \epsilon))| - | T \cap \Gamma(r(1-\epsilon))| $; we see this by noting that the dilation with factor $\frac{1 + \epsilon}{1 - \epsilon}$ mapping $\mathcal{F}(r(1 - \epsilon)) $ to $\mathcal{F}(r(1 + \epsilon))$ sweeps through $\mathcal{F}(r)$, and so $ \mathcal{F}(r) \subset \Gamma(r( 1 + \epsilon)) \backslash \Gamma(r( 1 - \epsilon)) $. It is therefore clear that

\begin{equation} \label{PlusMinusBound}
\lim_{r\to\infty} \frac{ |T \cap \mathcal{F}(r)|  }{ N(r) } \leq \lim_{r\to\infty} \left[\frac{ m(\Gamma) }{  m(B_n) }\right]\left[\frac{N(r(1 + \epsilon)) -  N(r(1 - \epsilon))}{N(r)}\right].
\end{equation}

\noindent However, the right-hand-side vanishes as $\epsilon \to 0$, so we are done. Thus, the contribution of points from each cube's boundary is negligible, and we henceforth neglect them in our considerations.   
\end{proof} 

\end{lemma}

Note that this easily implies that the contributions of points of $T \times T$ on finite $(2n - 1)$-dimensional hypersurfaces (such as the faces of $2n$-cubes) can be neglected as well. We will make use of both the $n$-dimensional and the $2n$-dimensional implications of the boundary points lemma in what follows.

\begin{lemma} 

Consider a $n$-box $C \subset \mathbb{R}^n$ whose sides are parallel to the axes and whose faces are closed, and let $C(r)$ be $C$ dilated by $r$ about $\mathbf{0}$. Then, we have that

\begin{equation} \lim_{r \to \infty} \frac{ \left| T \cap C(r) \right| }{N(r)} = \frac{m(C)}{m(B_n)}.
\label{CubeLemma}
\end{equation}

\begin{proof} 

Note that by the previous theorem, for any $\Omega$ with $\Vx \in \Omega \Longrightarrow \sigma \Vx \in \Omega \ \forall \ 0 \leq \sigma \leq 1$,
\[  \lim_{r \to \infty} \frac{\left| T \cap \Omega(r) \right|}{N(r)} = \frac{m(\Omega)}{m(B_{n})}.\]
Assume $C(r)$ does not contain the origin, for otherwise, we would be done by Theorem 2.2. Moreover, assume that no face of $C(r)$ passes through an axis of $\mathbb{R}^n$ (if one did, simply elongate $C(r)$ downward in that direction until the origin is reached. Then $C(r)$ will be the direct difference of the two boxes which pass through the origin, both of whose edges extend those of $C(r)$, and one of whose faces perpendicular to the axis in question is closer to the origin than the other). 

\vspace{2mm}

Otherwise, let us without loss of generality assume that $C(r)$ lies in the first quadrant. Now, let $\Vv_1$ be the vertex of $C(r)$ closest to the origin (for some large $r$), and let $\Vv_2$ be the vertex opposite $\Vv_1$ on the long diagonal. Define $U(r)$ to be the closed $n$-box with diagonally opposite vertices $\mathbf{0}$ and $\Vv_2$ and edges parallel to the coordinate axes, and let $\tilde{C}(r)$ be the closure of $U(r) \backslash C(r)$ (by Lemma 1, the number of points in $T$ on $\partial C(r)$ is negligible, so we may simply say $| T \cap \tilde{C}(r) | \sim | T \cap U(r) \backslash C(r) | $). Consider any point $\Vx$ on the boundary of $\tilde{C}(r)$. Then, it is easy to see that $\sigma\Vx \in \tilde{C}(r) \ \forall\ 0 \leq \sigma \leq 1$, since $\tilde{C}(r)$ extends all the way back to the origin in every direction $\mathbf{u}$ such that the components of $\Vu$ are nonnegative (and obviously every coordinate on the boundary of $\tilde{C}(r)$ is nonnegative). Therefore, we may apply Theorem 2.2 to $\tilde{C}(r)$. Note that we may obviously apply Theorem 2.2 to $U(r)$. We obtain

\begin{equation} \lim_{r \to \infty} \frac{\left| T \cap U(r) \right|}{N(r)} =  \lim_{r \to \infty} \frac{m(U(r))}{m(B_{n}(r))} 
\label{EquidistributionOnU}
\end{equation}
\begin{equation}
\lim_{r \to \infty} \frac{| T \cap \tilde{C}(r) |}{N(r)} = \lim_{r \to \infty} \frac{m(\tilde{C}(r) )}{m(B_{n}(r))} ; 
\label{EquidistributionOnCTilde}
\end{equation}
thus, because $| T \cap U(r) | - | T \cap U(r) \backslash C(r) | = |T \cap C(r)|$,
\[ \lim_{r \to \infty} \frac{\left| T \cap C(r) \right|}{N(r)} =  \lim_{r \to \infty} \frac{m(U(r)) - m(\tilde{C}(r) )}{m(B_{n}(r))} =  \lim_{r \to \infty} \frac{m(C(r))}{m(B_n(r))} = \frac{m(C)}{m(B_n)}. \]
\end{proof}

\end{lemma}

Now, dilate $\Psi$ and $\mathcal{H}_\delta$ by $r$ about $\mathbf{0}$; we consider a $2n$-cube $Q(r) \in \mathfrak{C}_{\delta}(r)$, where $\mathfrak{C}_{\delta}(r)$ consists of the cubes in $\mathfrak{C}_{\delta}$ dilated by $r$ about $\mathbf{0}$. Now, $Q(r) = Q_1(r) \times Q_2(r)$, where $Q_1(r)$ and $Q_2(r)$ are $n$-cubes in the space of $(x_1, \dots x_n)$ and $(x_{n+1}, \dots x_{2n})$, respectively. Consider a point in $Q_1(r)$ with $(x_1, \dots x_n) \in T$. Then, there are $\left| Q_2(r) \cap T \right|$ corresponding points in $Q_2(r)$ such that the entire point $(x_1, \dots x_{2n}) \in T \times T$. But there are $|Q_1(r) \cap T|$ such points in $Q_1(r)$ with $(x_1, \dots x_n) \in T$. Thus, we get $|Q_1(r) \cap T| |Q_2(r) \cap T| = | Q(r) \cap (T \times T) | $. Since $m(Q_1(r))m(Q_2(r)) = m(Q(r))$ (because all are Cartesian products of intervals), Lemma 2 gives that

\begin{equation} \label{EachBox} \lim_{r \to \infty} \frac{  | (T \times T) \cap Q(r)  | }{ N(r)^2} = \frac{m(Q)}{m(B_n)^2} .
\end{equation}
 
\noindent where we have used Lemma 1 in $2n$-dimensions to ignore the fact that some faces are open whereas others are closed. Summing over the $Q(r)$ in (\ref{EachBox}), we obtain
 
\begin{equation} \label{SumBox}
\lim_{r \to \infty} \frac{1}{N(r)^2} \left| (T \times T) \cap \Bigg(\bigcup_{Q\in\mathfrak{C}_{\delta}} Q(r)\Bigg)  \right| =  \sum_{Q \in \mathfrak{C}_{\delta}}\frac{m(Q)}{m(B_n)^2}.
\end{equation}

\noindent since the number of points omitted from consideration on the boundary of any of the $Q$ is negligible. Taking the limit as $\delta\to 0$, the union of all of the $Q$-regions in $\mathfrak{C}_{\delta}$ tends towards $\Psi$, so we have that

\begin{equation} \label{TcrossTResult} \lim_{r \to \infty}\frac{ | (T \times T) \cap \Psi(r) | }{ N(r)^2} = \frac{m(\Psi)}{m(B_n)^2} . 
\end{equation}

\noindent To finish, we simply note that, for each point $(x_1, \dots x_n)$ in $T\cap B_n(r)$, there are $N(r)$ points $(x_{n+1}, \dots x_{2n}) \in T$ such that $(x_1, \dots x_{2n}) \in K(r)$. There are $N(r)$ such points $(x_1, \dots x_n)$, so $|(T \times T) \cap K(r) |  = N(r)^2$, finishing the proof. 
\end{proof}
\end{theorem}

\section{Calculating the Pair Correlation Distribution}

Given the volume heuristic, we have all the mathematical framework we need to compute the pair correlation function for equidistributed sets. All that follows are a few volume calculations in $2n$-dimensions by means of which we can express the required distribution in terms of the incomplete beta functions (see below).

\subsection{Definition of Special Functions}

Before we begin calculating the pair correlation distribution, we first need to define some familiar special functions, which we will use to both calculate and represent the pair correlation distribution:

\subsubsection{Gamma Function}

The \textit{Gamma function} $\Gamma(z)$ is defined (\cite[Eq.~5.2.1]{DLMF}) as

\begin{equation} \label{GammaDefinition}
\Gamma(z) \equiv \int_0^{\infty} e^{-t} t^{z-1}\ dt.
\end{equation}

\noindent It takes the value $(n-1)!$ at the positive integer $n$, and the value $\sqrt{\pi}$ at $1/2$ (\cite[Eq.~5.4.1]{DLMF},\cite[Eq.~5.4.6]{DLMF}). Furthermore, it obeys the following duplication formula:

\begin{equation} \label{GammaDuplication}
\Gamma(2z)=\frac{1}{\sqrt{\pi}} 2^{2z-1} \Gamma(z)\Gamma\left(z+\frac{1}{2}\right).
\end{equation}

\subsubsection{Beta Function}

The \textit{Beta function} $\Beta(a,b)$ is defined (\cite[Eq.~5.12.1]{DLMF}) as 

\begin{equation} \label{BetaDefinition}
\Beta(a,b) \equiv \int_0^1 t^{a-1} (1-t)^{b-1}\ dt = \frac{\Gamma(a)\Gamma(b)}{\Gamma(a+b)}.
\end{equation}

\noindent It is also equal to (\cite[Eq.~5.12.2]{DLMF})

\begin{equation} \label{BetaTrig}
\Beta(a,b) = 2\int_0^{\frac{\pi}{2}} \sin^{2a-1} \theta \cos^{2b-1}\theta\ d\theta.
\end{equation}

\subsubsection{Incomplete Beta Functions}

The \textit{incomplete Beta function} $\Beta_x(a,b)$ is defined as (\cite[Eq.~8.17.1]{DLMF})

\begin{equation} \label{IncompleteBetaDefinition}
\Beta_x(a,b) \equiv \int_0^x t^{a-1} (1-t)^{b-1}\ dt.
\end{equation}

\noindent It is also useful to define the \textit{regularized incomplete Beta function} (\cite[Eq.~8.17.2]{DLMF})

\begin{equation} \label{RegularizedIncompleteBetaDefinition}
I_x(a,b) \equiv \frac{\Beta_x(a,b)}{\Beta(a,b)}.
\end{equation}

\noindent This obeys the reflection formula (\cite[Eq.~8.17.4]{DLMF})

\begin{equation} \label{BetaReflection}
I_x(a,b)+I_{1-x}(b,a)=1.
\end{equation}

\subsection{Evaluation of the Volume of the Region}

From the discussion in the introduction, the problem of pair correlations reduces to the following: We want to find the number of pairs of points $(a_1, a_2, ..., a_n)$ and $(b_1, b_2, ..., b_n)$ in  $T$ satisfying 

\begin{equation} \label{MagnitudeConditions}
\sum_{i=1}^n a_i^2 \leq R^2,\ \ \sum_{i=1}^n b_i^2 \leq R^2 ,
\end{equation}

\begin{equation} \label{DistanceConditions}
\sum_{i=1}^n \left(a_i-b_i\right)^2 \leq \lambda^2 R^2 
\end{equation}

\noindent for given $\lambda, R$. These conditions can be visualized as the intersection of three regions in $\mathbb{R}^{2n}$. Now, we are in a position to apply Theorem 2.3: the number of ordered pairs of points in $T$ within the required compact region $\Psi(R, \lambda)$ (defined by the above two equations) can be asymptotically $(R \to \infty)$ estimated by the following: 

\[ | (T \times T) \cap \Psi(R, \lambda) | \sim  N(R)^2 \frac{  m (\Psi(R, \lambda)) }{m(B_n)^2}  . \]

\noindent However, note that $\Psi(R, 2) = B_n(R) \times B_n(R)$. Thus, 

\[ \frac{ | (T \times T) \cap \Psi(R, \lambda) | }{ | (T \times T) \cap \Psi(R, 2) | } \sim \frac{  m (\Psi(R, \lambda)) }{m(\Psi(R, 2)) } . \]

Therefore, the calculated distribution function of $\lambda$ will be the same for both $| (T \times T) \cap \Psi(R, \lambda) |$ and $m (\Psi(R, \lambda)) $. 
We therefore completely disregard $T$ and focus only on computing the volume of $\Psi(R, \lambda)$. When computing a probability distribution, the scaling radius $R$ becomes irrelevant, so we without loss of generality we let $R = 1$. 

We make the change of variables

\begin{equation} \label{ChangeOfVars}
a_i = v_i+\frac{u_i}{2},\ \ b_i = v_i-\frac{u_i}{2}.
\end{equation}

\noindent It is easy to verify that the Jacobian of this map is $1$. Note that $a_i-b_i = u_i$, so condition (\ref{DistanceConditions}) becomes

\[\sum_{i=1}^n u_i^2 \leq \lambda^2\ \Leftrightarrow\ |\Vu|\leq \lambda,\]

\noindent and the conditions in (\ref{MagnitudeConditions}) can be rephrased as requiring that the point $\mathbf{v} = (v_1, ..., v_n)$ is inside both of the $n$-balls of radius $1$ centered at each of 

\begin{equation} \label{ChangeOfVarsNewCenters}
\mathbf{-u} = \left(-\frac{u_1}{2},...,-\frac{u_n}{2}\right)\ \mathrm{and}\ \mathbf{u} = \left(\frac{u_1}{2},...,\frac{u_n}{2}\right).
\end{equation}

\noindent The volume of this region is just twice the volume of the hyperspherical cap with radius $1$ and cap base height $\frac{|\mathbf{u}|}{2}$. (We define $r = |\mathbf{u}|$.) By the theorem in \cite{Li2011}, this volume is just

\begin{equation} \label{CapVolume}
\frac{\pi^{\frac{n}{2}}}{2\Gamma\left(\frac{n}{2}+1\right)}I_{1-\frac{r^2}{4}}\left(\frac{n+1}{2}, \frac{1}{2}\right).
\end{equation}

\noindent Thus, the volume of the region defined by a given vector $\mathbf{u}$ is 

\begin{equation} \label{CapVolumeIntegral}
\frac{\pi^{\frac{n}{2}}}{\Gamma\left(\frac{n+1}{2}\right)\Gamma\left(\frac{1}{2}\right)}\int_0^{1-\frac{r^2}{4}} t^{\frac{n-1}{2}}(1-t)^{-\frac{1}{2}} dt.
\end{equation}

\noindent To obtain the volume of the region described by (\ref{MagnitudeConditions}) and (\ref{DistanceConditions}), we must integrate this over the $n$-ball with radius $\lambda$. However, since our function depends only on $r$, we can multiply by the surface area of each hyper-spherical shell of radius $r$, and integrate from $r=0$ to $\lambda$. This surface area formula is well-known (and found in \cite{Li2011}) to be

\begin{equation} \label{SurfaceArea}
\frac{2\pi^{\frac{n}{2}}}{\Gamma\left(\frac{n}{2}\right)}r^{n-1},
\end{equation}

\noindent and thus the volume is

\[\frac{2\pi^n}{\Gamma\left(\frac{n}{2}\right)\Gamma\left(\frac{n+1}{2}\right)\Gamma\left(\frac{1}{2}\right)}\int_0^\lambda r^{n-1}\int_0^{1-\frac{r^2}{4}} t^{\frac{n-1}{2}}(1-t)^{-\frac{1}{2}} dt\ dr.\]

\noindent Interchanging the order of integration and simplifying using the incomplete Beta function, the expression of the volume reduces to the following:
 
\begin{equation} \label{LambdaVolume}
\frac{2\pi^n}{n\Gamma\left(\frac{n}{2}\right)\Gamma\left(\frac{n+1}{2}\right)\Gamma\left(\frac{1}{2}\right)}
\Bigg[ \lambda^n\Beta_{1-\frac{\lambda^2}{4}}\left( \frac{n+1}{2}, \frac{1}{2}\right) \Bigg]  +  
\end{equation}

\[ + \frac{2\pi^n}{n\Gamma\left(\frac{n}{2}\right)\Gamma\left(\frac{n+1}{2}\right)\Gamma\left(\frac{1}{2}\right)} \Bigg[ 2^n\Beta_{\frac{\lambda^2}{4}}\left( \frac{n+1}{2}, \frac{n+1}{2}\right) \Bigg] . \]

\subsection{Deriving a Probability Density Function}

To find the constant by which we need to divide to turn (\ref{LambdaVolume}) into a proper cumulative density function, we evaluate (\ref{LambdaVolume}) at $\lambda = 2$, which measures the volume of the whole region (given by (\ref{MagnitudeConditions}) only):
 
\begin{equation} \label{EntireVolume}
\frac{2\pi^n}{n\Gamma\left(\frac{n}{2}\right)\Gamma\left(\frac{n+1}{2}\right)\Gamma\left(\frac{1}{2}\right)}
2^n\Beta\left(1; \frac{n+1}{2}, \frac{n+1}{2}\right) = \frac{\pi^n2^n\Gamma\left(\frac{n+1}{2}\right)}{\Gamma\left(\frac{n}{2}+1\right)\sqrt{\pi}\Gamma(n+1)}. 
\end{equation}

\noindent The duplication formula (\ref{GammaDuplication}) simplifies this to 

\[\frac{\pi^n}{\Gamma\left(\frac{n}{2}+1\right)^2}, 
\]

\noindent which is the square of the volume of an $n$-ball of radius $1$ (found in \cite{Li2011}), as expected (since each of the normalized $\mathbf{a}$ and $\mathbf{b}$ vectors are in the space of the $n$-ball of radius $1$). Dividing the expression (\ref{LambdaVolume}) by this, we get

\begin{equation} \label{LambdaProbability}
\frac{1}{\Beta\left(\frac{n+1}{2}, \frac{1}{2}\right)}
\left(\lambda^n\Beta\left(1-\frac{\lambda^2}{4}; \frac{n+1}{2}, \frac{1}{2}\right) + 2^n\Beta\left(\frac{\lambda^2}{4}; \frac{n+1}{2}, \frac{n+1}{2}\right)\right),
\end{equation}

\noindent which is the probability that the distance between any two given lattice points in the $n$-ball with radius $R$ (sufficiently large) is less than $R\lambda$. To obtain the probability density function, we differentiate with respect to $\lambda$:

\begin{equation} \label{FinalProbability}
P_n(\lambda) = n\lambda^{n-1}I_{1-\frac{\lambda^2}{4}}\left(\frac{n+1}{2},\frac{1}{2}\right).
\end{equation}

\subsection{Specific Values of $n$}

\noindent For $n=2$, the function simplifies to

\begin{equation} \label{n=2}
\frac{4\lambda}{\pi}\cos^{-1}\left(\frac{\lambda}{2}\right) - \frac{\lambda^2\sqrt{4-\lambda^2}}{\pi},
\end{equation}

\noindent as has been very well documented. For $n=3$, the function simplifies to 

\begin{equation} \label{n=3}
3\lambda^2-\frac{9\lambda^3}{4}+\frac{3\lambda^5}{16}.
\end{equation}

\noindent Note that, for odd $n$, the probability density function will be a polynomial in $\lambda$ of degree $2n-1$ with rational coefficients.
 
\vspace{2mm}

\noindent We now examine what happens to the probability density function as $n$ tends to $\infty$. Consider the limit of the integral (\ref{FinalProbability}) as $n \to \infty$, using the substitution $u=(1-t)^{-1/2}$:

\begin{equation} \label{LimitIntegralStatement}
P_{\infty}(\lambda) \equiv \lim_{n\to\infty} \frac{2n\lambda^{n-1}}{\Beta\left(\frac{n+1}{2},\frac{1}{2}\right)}\int_{\frac{\lambda}{2}}^1 \left(1-u^2\right)^{\frac{n-1}{2}}\ du.
\end{equation}

\noindent We use the method of Laplace in \cite{Laplace1820} for asymptotically estimating integrals to obtain

\begin{equation} \label{LimitSubst1}
P_{\infty} (\lambda) = \lim_{n\to\infty} C(\lambda)\frac{n\lambda^{n-1}}{\Beta\left(\frac{n+1}{2},\frac{1}{2}\right)}\sqrt{\frac{2\pi}{\frac{n-1}{2}}}\left(1-\frac{\lambda^2}{4}\right)^{\frac{n-1}{2}},
\end{equation}

\noindent where $C(\lambda)$ is some positive value independent of $n$. By Stirling's Approximation, we have that

\begin{equation} \label{PinftyAsymptotics}
P_{\infty} (\lambda) = C(\lambda)\sqrt{2}\lim_{n\to\infty} n\lambda^{n-1}\left(1-\frac{\lambda^2}{4}\right)^{\frac{n-1}{2}}.
\end{equation}

\noindent Letting $\lambda=2\sin\theta$ simplifies this to $C(\lambda)\sqrt{2}\ n\sin^{n-1}(2\theta)$, which tends to $0$ everywhere except $\lambda=\sqrt{2}$, where it tends to $\infty$. Thus, this probability density function can be said to equal $\delta(\lambda-\sqrt{2})$, where $\delta$ is the Dirac delta function.

\section{Specific sets $T\subset \mathbb{R}^n$}

In \cite{Wills1973}, \cite{Steinhaus1947}, and \cite{Clark2013}, it is proven that for any subset $\Omega\subset \mathbb{R}^n$, the number of lattice points or primitive lattice points inside $\Omega$ dilated by radius $r$ grows like $R^n$ times a constant ($1$ for lattice points and $1/\zeta(n)$ for primitive lattice points). Therefore, by Theorem 2.1, the probability density function (\ref{FinalProbability}) is also the probability density function for the normalized distance between pairs of lattice and primitive lattice points.

\subsection{Numerical Results}
We may verify with simple numerical checks that the probability density (\ref{FinalProbability}) agrees extremely well with the actual distribution of normalized distances between lattice and primitive lattices points enclosed within an $n$-ball of radius $R$. Here, we present results for $n = 2$, $R = 300$ and $n = 3$, $R = 30$: 

\begin{figure}[H]

\centering

\includegraphics[scale=1.20]{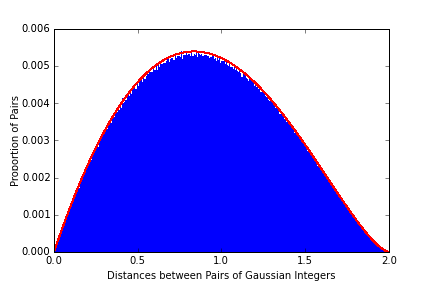} 

\caption{This graph shows the distribution of distances for $n = 2$ and $R = 300$ over the lattice points (Gaussian integers). The histogram (blue) plots the relative frequency with which a certain normalized distance occurs. The curve (red) represents the (appropriately scaled) probability distribution given by (\ref{FinalProbability}). }

\end{figure}

\begin{figure}[H]

\centering

\includegraphics[scale=1.20]{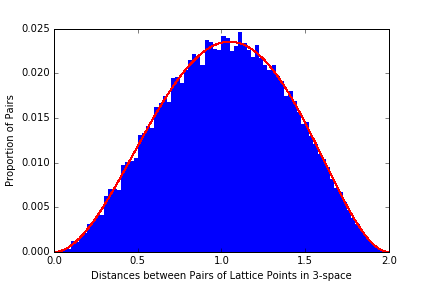} 

\caption{This graph shows the distribution of distances for $n = 3$ and $R = 30$ over the lattice points. }

\end{figure}

\begin{figure}[H]

\centering

\includegraphics[scale = 1.20]{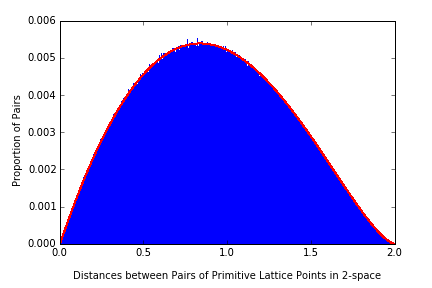}

\caption{This graph shows the distribution of distances for $n = 2$ and $R = 300$ over the primitive lattice points. }

\end{figure}

\begin{figure}[H]

\centering

\includegraphics[scale = 1.20]{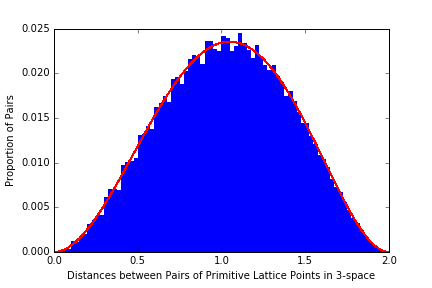}

\caption{This graph shows the distribution of distances for $n = 3$ and $R = 30$ over the primitive lattice points. }

\end{figure}

\section*{Acknowledgements}

We would like to thank Dr. Jayadev Athreya at the University of Washington for his continual support and encouragement throughout the course of our research.

\bibliography{PairCorrelations}

\end{document}